\documentclass[reqno]{amsart}
\usepackage{amssymb,eucal}

\newcommand{\jz}{\ensuremath{\set{\vee,0}}}
\newcommand{\one}{\mathbf{1}}
\newcommand{\two}{\mathbf{2}}

\newcommand{\jirr}{join-ir\-re\-duc\-i\-ble}

\newcommand{\mirr}{meet-ir\-re\-duc\-i\-ble}
\newcommand{\Mirr}{Meet-ir\-re\-duc\-i\-ble}
\newcommand{\jsd}{join-sem\-i\-dis\-trib\-u\-tive}

\newcommand{\fbs}{finitely bi-spatial}
\newcommand{\Fbs}{Finitely bi-spatial}

\newcommand{\pup}[1]{\textup{(}#1\textup{)}}

\newcommand{\set}[1]{\{#1\}}
\newcommand{\setm}[2]{\{#1\mid#2\}}

\newcommand{\es}{\varnothing}

\newcommand{\id}{\mathrm{id}}

\newcommand{\ol}[1]{\overline{#1}}

\DeclareMathOperator{\Neu}{Neu}
\DeclareMathOperator{\Cen}{Cen}

\newcommand{\SP}{\mathbf{S_p}}
\newcommand{\Co}{\mathbf{Co}}

\DeclareMathOperator{\SPl}{S_p}

\DeclareMathOperator{\J}{J}
\DeclareMathOperator{\M}{M}

\newcommand{\Pow}{\mathcal{P}}
\newcommand{\PO}{\mathcal{P}(\omega)}
\newcommand{\POF}{\mathcal{P}(\omega)/{\mathrm{fin}}}

\theoremstyle{plain}

\newtheorem{lemma}{Lemma}[section]
\newtheorem{theorem}[lemma]{Theorem}
\newtheorem{proposition}[lemma]{Proposition}
\newtheorem{corollary}[lemma]{Corollary}
\newtheorem{example}[lemma]{Example}

\newtheorem*{stat}{\name}
\newcommand{\name}{testing}

\theoremstyle{definition}
\newtheorem{definition}[lemma]{Definition}
\newtheorem*{problem}{Problem}
\newtheorem*{notation}{Notation}

\theoremstyle{remark}
\newtheorem{remark}[lemma]{Remark}

\newcommand{\qedc}{{\qed}~{\rm Claim~{\theclaim}.}}

\begin{document}

\title[Direct decompositions of lattices]
{Direct decompositions of non-algebraic\\
complete lattices}

\author[F.~Wehrung]{Friedrich Wehrung}
\address{CNRS, UMR 6139\\
D\'epartement de Math\'ematiques\\
Universit\'e de Caen\\
14032 Caen Cedex\\
France}
\email{wehrung@math.unicaen.fr}
\urladdr{http://www.math.unicaen.fr/\~{}wehrung}


\date{\today}

\begin{abstract}
For a given complete lattice $L$, we investigate whether $L$ can be
decomposed as a direct product of directly indecomposable lattices. We prove
that this is the case if every element of $L$ is a join of \jirr\ 
elements and dually, thus extending to non-algebraic lattices a result of
L.~Libkin. We illustrate this by various examples and counterexamples.
\end{abstract}

\maketitle

\section{Introduction}
L. Libkin proves in \cite{Libk95} that if an algebraic lattice $L$ is
\emph{spatial}, that is, every element of $L$ is a join of completely
join-irreducible elements of $L$, then $L$ can be decomposed as a direct
product of directly indecomposable lattices---we say that $L$ is
\emph{totally decomposable}. This result extends the classical one about
decomposing a geometric lattice as a product of indecomposable factors. It
is in turn extended in J.~Jakubik~\cite{Jaku96} by relaxing the
completeness assumptions on $L$, and in A.~Walendziak \cite{Wale97} to
algebraic lattices in which the unit element is a join of \jirr\ elements.
None of these results avoids the assumption that the lattice is compactly
generated, in particular, they do not apply to the closure
lattices of the so-called
\emph{convex geometries} studied in~\cite{AGT}, as the latter are not
algebraic as a rule (by definition, a convex geometry is a
closure space satisfying the anti-exchange property).

In this paper, we extend Libkin's methods and result to a class of
lattices that properly contains both Libkin's lattices and all closure
lattices of most convex geometries, the class of \emph{\fbs} complete
lattices (Definition~\ref{D:finbisp}), see Theorem~\ref{T:DecompL}. 
We also illustrate this by a few examples and counterexamples that show,
in particular, that our assumptions cannot be relaxed much:
\begin{itemize}\em
\item There exists a self-dual, complete, distributive lattice $D$ whose
center is a complete atomistic sublattice but $D$ is not totally
decomposable \pup{see Example~\textup{\ref{Ex:NonDec}}}.

\item There exists a dually algebraic, atomistic, distributive lattice
whose center is not complete \pup{see Example~\textup{\ref{Ex:DCX}}}.

\item Denote by $\SP(A)$ the lattice of algebraic subsets of a complete
lattice $A$. If $A$ is Boolean, then $\SP(A)$ is subdirectly irreducible
\pup{Proposition~\textup{\ref{P:SPIndec}}}, but for $A$ a chain, $\SP(A)$
may not have complete center
\pup{see Example~\textup{\ref{Ex:NonCpleteCen}}}.
\end{itemize}

We observe that Examples~\ref{Ex:DCX} and \ref{Ex:NonCpleteCen} solve
negatively a problem formulated by M.\,F. Janowitz in \cite{Jano67},
whether the center of a complete lattice must be a complete sublattice.

For a set $X$, we denote by $\Pow(X)$ the powerset lattice of $X$. We adopt
the standard set-theoretical notation for ordinals, for example,
$n=\set{0,\ldots,n-1}$ for every nonnegative integer $n$,
then $\omega=\set{0,1,2,\ldots}$, and $\omega+1=\omega\cup\set{\omega}$.

An element $p$ of a lattice $L$ is \emph{\jirr} (resp., \emph{completely
\jirr}), if it is nonzero if $L$ has a zero, and $p=x\vee y$ implies that
$p\in\set{x,y}$, for all $x$, $y\in L$ (resp., $p$ has a unique lower
cover). \Mirr\ (resp., completely \mirr) elements are defined dually.
We denote by $\J(L)$ (resp., $\M(L)$) the set of all \jirr\ (resp., \mirr)
elements of a lattice $L$.

For elements $x$ and $y$ of a given poset, let $x\prec y$ be the statement
that $x<y$ and there is no element strictly between $x$ and $y$.
A lattice $L$ with zero is \emph{atomistic}, if every element of $L$ is a
join of atoms of $L$.

\section{Decompositions of complete lattices}\label{S:DecCpl}

We first recall some standard terminology and facts, see
\cite[Chapter~III, Section~2]{GLT2}.
An element $a$ in a lattice $L$ is \emph{neutral}, if $\set{a,x,y}$
generates a distributive sublattice of $L$, for all $x$, $y\in L$. We shall
denote by $\Neu L$ the subset of all neutral elements of $L$. If $L$ is
\emph{bounded}, we say that an element $a$ of $L$ is \emph{central}, if it is
both neutral and complemented in $L$; then the complement $\neg a$ is
unique, and it is also central. Hence $\Neu L$ is a distributive sublattice
of~$L$, and, if $L$ is bounded, then $\Cen L$ is a Boolean sublattice of
$L$.

The elements of $\Cen L$ correspond exactly to the direct decompositions of
$L$. This can be expressed conveniently in the following way, see
\cite[Theorem~III.4.1]{GLT2}:

\begin{lemma}\label{L:charactC(L)}
Let $L$ be a bounded lattice, let $a$, $b\in L$. Then the following are
equivalent:
\begin{enumerate}
\item There are bounded lattices $A$ and $B$ and an isomorphism
$f\colon L\to A\times B$ such that $f(a)=(1,0)$ and $f(b)=(0,1)$.

\item $(a,b)$ is a complementary pair of elements of $\Cen L$, that is,
$a$, $b\in\Cen L$, $a\nobreak\wedge\nobreak b\nobreak=\nobreak0$, and
$a\vee b=1$.
\end{enumerate}
\end{lemma}

We observe the following easy consequence of Lemma~\ref{L:charactC(L)}:

\begin{proposition}\label{P:At(C(L))}
Let $L$ be a bounded lattice, let $a\in\Cen L$. Then the following assertions
hold:
\begin{enumerate}
\item $\Cen([0,a])=\Cen L\cap[0,a]$.

\item If $a$ is an atom of $\Cen L$, then the interval $[0,a]$ is directly
indecomposable.
\end{enumerate}
\end{proposition}

\begin{definition}
A lattice $L$ is \emph{totally decomposable}, if it is isomorphic to a direct
product of the form $\prod_{i\in I}L_i$, where all the $L_i$-s are directly
indecomposable.
\end{definition}

Totally decomposable complete lattices can be easily characterized as follows:

\begin{proposition}\label{P:TotDec}
Let $L$ be a complete lattice. Then the following are equivalent:
\begin{enumerate}
\item $L$ is totally decomposable;

\item $\Cen L$ is a complete sublattice of $L$, it is atomistic, and, if $U$
denotes the set of its atoms, then the following holds:
 \begin{equation}\label{Eq:Decompx}
 x=\bigvee_{u\in U}(x\wedge u),\quad\text{for all }x\in L.\tag{J}
 \end{equation}
\end{enumerate}
\end{proposition}

\begin{proof}
(i)$\Rightarrow$(ii) Suppose that $L=\prod_{i\in I}L_i$, for a family
$(L_i)_{i\in I}$ of directly indecomposable lattices. Observe that all the
$L_i$-s are complete, in particular, they are bounded lattices. For all
$X\subseteq I$, the characteristic function $\chi_X$ of $X$ in $I$ belongs to
the center of $L$, and its complement is $\chi_{I\setminus X}$. The
complemented pair $(\chi_X,\chi_{I\setminus X})$ of elements of $\Cen L$
induces an isomorphism $L\cong L_X\times L_{I\setminus X}$, where we put
$L_Y=\prod_{i\in Y}L_i$ for every subset $Y$ of $I$. Conversely, if
$u=(u_i)_{i\in I}$ is an element of $\Cen L$, then
$u_i\in\Cen L_i$, for all $i\in I$, thus, since $L_i$ is directly
indecomposable, $u_i\in\set{0,1}$. Therefore, $u=\chi_X$, where
$X=\setm{i\in I}{u_i=1}$.

Consequently, $\Cen L=\setm{\chi_X}{X\subseteq I}$ is a complete sublattice of
$L$. Furthermore, it is atomistic, with atoms the elements $\chi_{\set{i}}$
for $i\in I$. The assertion \eqref{Eq:Decompx} follows easily.

(ii)$\Rightarrow$(i) Suppose that (ii) holds, and denote by $U$ the set of
all atoms of $\Cen L$. Put $L_u=[0,u]$, for all $u\in U$,
then $L'=\prod_{u\in U}L_u$, and define maps
$f\colon L\to L'$ and $g\colon L'\to L$ by the
rules
 \begin{align*}
 f(x)&=(x\wedge u)_{u\in U},&&\text{for all }x\in L,\\
 g\bigl((x_u)_{u\in U}\bigr)&=\bigvee_{u\in U}x_u,&&
 \text{for all }(x_u)_{u\in U}\in L'.
 \end{align*}
For $(x_u)_{u\in U}\in L'$, if we put $x=\bigvee_{u\in U}x_u$,
then, for any $u\in U$, we obtain, by using the fact that $u$ is neutral, the
inequalities $x_u\leq x\wedge u\leq(x_u\vee\neg u)\wedge u=x_u$, whence
$x_u=x\wedge u$. Hence $f\circ g=\id_{L'}$.
Moreover, $g\circ f=\id_L$ follows from the assumption~(ii). Hence, $f$ and
$g$ are mutually inverse isomorphisms. By Proposition~\ref{P:At(C(L))}(ii),
all the factors of the form $L_u$ are directly indecomposable.
\end{proof}

\begin{remark}\label{R:TotDec}
For a bounded lattice $L$, the completeness assumption in
Proposition~\ref{P:TotDec} can be much relaxed. For example,
Proposition~\ref{P:TotDec} remains valid under the assumption that any
family $(x_u)_{u\in U}$ with $x_u\leq u$, for all $u\in U$, has a
join, and the proof is the same.
\end{remark}

In our next result, we shall state a number of conditions that imply
\eqref{Eq:Decompx}. In order to state it conveniently, we set a definition,
that will also be used in Section~\ref{S:Decomp}:

\begin{definition}\label{D:spatial}
Let $L$ be a lattice. We say that $L$ is \emph{finitely spatial} (resp.,
\emph{spatial}), if every element of $L$ is a join of \jirr\ (resp.,
completely \jirr) elements of~$L$. Let \emph{dually spatial}, resp.
\emph{dually finitely spatial}, be the dual notions.
\end{definition}

For example, the real unit interval $[0,1]$ is finitely spatial but not
spatial.

\begin{proposition}\label{P:SuffDec}
Let $L$ be a complete lattice such that $\Cen L$ is a complete atomistic
sublattice of $L$. Then each of the following conditions \pup{and also its
dual} implies that $L$ is totally decomposable:
\begin{enumerate}
\item $L$ is upper continuous.

\item $L$ is \emph{separative}, that is, for any elements $x$, $y\in L$
such that $x\nleq y$, there exists $z\in L$ such that $0<z\leq x$ and
$z\wedge y=0$.

\item $L$ is finitely spatial.
\end{enumerate}
\end{proposition}

\begin{proof}
By Proposition~\ref{P:TotDec}, it suffices that the condition
\eqref{Eq:Decompx} is satisfied by $L$. So let $x\in L$. We put
$y=\bigvee_{u\in U}(x\wedge u)$. In case $L$ is upper continuous, we observe
that $x\wedge\bigvee V=\bigvee_{u\in V}(x\wedge u)$ for every finite subset
$V$ of $U$ (because all elements of $U$ are neutral). Hence, by the upper
continuity of $L$,
 \[
 x=x\wedge\bigvee U=\bigvee_{V\subseteq U\,\text{finite}}
 \left(x\wedge\bigvee V\right)
 =\bigvee_{u\in U}(x\wedge u)=y.
 \]
We conclude the proof of (i) by Proposition~\ref{P:TotDec}.

Suppose that $L$ is separative and that $y<x$. Then, by
assumption, there exists $z\in L$ such that $0<z\leq x$ but $z\wedge y=0$.
Hence, for all $u\in U$, the equality $z\wedge u=0$ holds, thus $z\leq\neg u$.
Therefore, $z\leq\bigwedge_{u\in U}\neg u=\neg\bigvee U=0$, a contradiction.

Finally, suppose that $L$ is finitely spatial. To prove that
\eqref{Eq:Decompx} holds at all elements of $L$, it suffices to verify it for
$x\in\J(L)$. Suppose that it is not the case, that is,
$x>\bigvee_{u\in U}(x\wedge u)$, where $U$ denotes the set of atoms of
$\Cen L$. Every element $u$ of~$U$ belongs to $\Cen L$, whence
$x=(x\wedge u)\vee(x\wedge\neg u)$, but $x\wedge u<x$ by assumption and $x$
is \jirr, thus $x\wedge\neg u=x$, that is, $x\leq\neg u$. This holds for
all $u\in U$, therefore, by assumption on $\Cen L$, $x=0$, a contradiction.
\end{proof}

In particular, we observe that condition (ii) of Proposition~\ref{P:SuffDec}
holds if $L$ is either \emph{atomistic} or \emph{sectionally complemented}.
Since the center of a complete relatively complemented lattice is a complete
sublattice, see \cite{Jano67}, we obtain the following result:

\begin{corollary}\label{C:DecRelCpl}
Let $L$ be a complete relatively complemented lattice. If $\Cen L$ is
atomistic, then $L$ is totally decomposable.
\end{corollary}

To conclude the present section, we shall now see that the condition
\eqref{Eq:Decompx} is not redundant in the statement of
Proposition~\ref{P:TotDec}. 

\begin{example}\label{Ex:NonDec}
There exists a self-dual, complete, distributive lattice $D$ such that
$\Cen D$ is a complete atomistic sublattice of $D$ but $D$ is not totally
decomposable.
\end{example}

\begin{proof}
{}From the classical theory of Boolean algebras, we know that any Boolean
algebra can be embedded into a complete Boolean algebra, see, for example,
\cite[Lemma~II.4.12]{GLT2}. We apply this to the
Boolean algebra $\POF$ of all subsets of $\omega$ modulo the ideal of finite
subsets, to embed it into a complete Boolean algebra $B$. We denote by $[x]$
the equivalence class, modulo the ideal of finite sets, of any subset $x$ of
$\omega$. We observe that $x\mapsto[x]$ defines a homomorphism of Boolean
algebras from $\PO$ to $B$. Thus, the subset $D$ of $\PO\times B\times\PO$
defined as
 \[
 D=\setm{(x,\alpha,y)\in\PO\times B\times\PO}
 {x\subseteq y\text{ and }[x]\leq\alpha\leq[y]}
 \]
is a sublattice of $\PO\times B\times\PO$, in particular, it is a
distributive lattice. Furthermore, $D$ is self-dual, \emph{via} the map
$(x,\alpha,y)\mapsto(\omega\setminus y,\neg\alpha,\omega\setminus x)$.

Let $\varphi\colon\PO\to D$, $x\mapsto(x,[x],x)$. It is obvious that
$\varphi$ is a $0,1$-lattice embedding. Furthermore, since $D$ is a bounded
distributive lattice, the center of $D$ consists exactly of the complemented
elements of $D$. Since $\varphi$ is a $0,1$-lattice homomorphism from $\PO$
to $D$, the range of $\varphi$ is contained in the center of $D$. Conversely,
if $z=(x,\alpha,y)$ is an element of $\Cen D$, then $z$ has a complement,
say, $z'=(x',\alpha',y')\in D$, so $x'=\omega\setminus x$ and
$y'=\omega\setminus y$, thus, since $x\subseteq y$ and $x'\subseteq y'$, we
obtain that $x=y$ and $x'=y'$, whence $\alpha=[x]$, so $z=\varphi(x)$.
Therefore, $\Cen D$ is the range of $\varphi$. It is atomistic, with
atoms the elements $a_n=\varphi(\set{n})=(\set{n},0,\set{n})$, for
$n<\omega$.

We now claim that $D$ is a complete lattice.
Indeed, let $(x_i,\alpha_i,y_i)_{i\in I}$ be a family of elements of~$D$, we
prove that it has a greatest lower bound in $D$. Put $x=\bigcap_{i\in I}x_i$,
$y=\bigcap_{i\in I}y_i$, and $\alpha=\bigwedge_{i\in I}\alpha_i\wedge[y]$. It
is obvious that $(x,\alpha,y)$ belongs to $D$ and that it is contained in
$(x_i,\alpha_i,y_i)$, for all $i\in I$. Let $(x',\alpha',y')\in D$ such
that $(x',\alpha',y')\leq(x_i,\alpha_i,y_i)$, for all $i\in I$. Then
$x'\subseteq x$ and $y'\subseteq y$, thus, since $\alpha'\leq\alpha_i$,
for all $i\in I$, and $\alpha'\leq[y']\leq[y]$, we obtain that
$\alpha'\leq\alpha$. So we have verified that $(x,\alpha,y)$ is the
greatest lower bound of $\setm{(x_i,\alpha_i,y_i)}{i\in I}$ in $D$; whence
$D$ is a complete lattice.

Moreover, in the particular case where $x_i=y_i$, for all $i\in I$ (so
$\alpha_i=[x_i]$), we obtain that $(x,\alpha,y)=(x,[x],x)$, where
$x=\bigcap_{i\in I}x_i$. Hence, $\varphi$ is a complete meet embedding.
The verification of the fact that $\varphi$ is a complete join embedding is
similar. Hence, $\varphi$ is a complete lattice embedding from $\PO$ into
$D$. Therefore, the center of $D$, which is also the range of $\varphi$, is a
complete sublattice of $D$.

Now put $b=(\es,1,\omega)$ (so $b\in D$).
We observe that $b\wedge a_n=(\es,0,\set{n})$, for all $n<\omega$,
hence
 \[
 \bigvee_{n<\omega}(b\wedge a_n)=(\es,0,\omega)<b.
 \]
By Proposition~\ref{P:TotDec}, $D$ is not totally decomposable.
\end{proof}

\begin{remark}
It is easy to read, in the proof above, the places where
Example~\ref{Ex:NonDec} fails the conditions (i)--(iii) of
Proposition~\ref{P:SuffDec}. For all $n<\omega$, the element
$\ol{a}_n=(n,0,n)$ belongs to $D$, while $\bigvee_{n<\omega}\ol{a}_n=1$ and
$\bigvee_{n<\omega}(\ol{a}_n\wedge b)<b$, thus verifying that $D$ is not
upper continuous. Put $\ol{b}=(\es,0,\omega)$. Then $\ol{b}<b$, while there is no
nonzero $z\leq b$ such that $z\wedge\ol{b}=0$, thus verifying that $D$ is
not separative. Finally, the \jirr\ elements below $b$ are exactly all the
$(\es,0,\set{n})$, and these join to $\ol{b}<b$, thus verifying that $D$
is not finitely spatial.

\end{remark}

\section{\Fbs\ complete lattices}
\label{S:Decomp}

We start by defining the objects of the section title:

\begin{definition}\label{D:finbisp}
We say that a bounded lattice $L$ is \emph{\fbs}, if it is both finitely
spatial and dually finitely spatial (see Definition~\ref{D:spatial}).
\end{definition}

\begin{notation}
Let $x\in L$, let $(x_i)_{i\in I}$ be a family of elements
of $L$. Let $x=\bigvee^*_{i\in I}x_i$ hold, if
 \[
 p\leq x\quad\text{if{f}}\quad\exists i\in I\text{ such that }
 p\leq x_i,\text{ for all }p\in\J(L).
 \]
For $|I|=2$, we define similarly the notation $z=x\vee^*y$, for $x$, $y$,
$z\in L$. Similarly, let $x=\bigwedge^*_{i\in I}x_i$ hold, if
 \[
 x\leq u\quad\text{if{f}}\quad\exists i\in I\text{ such that }
 x_i\leq u,\text{ for all }u\in\M(L),
 \]
and, for $|I|=2$, we define similarly the notation $z=x\wedge^*y$.
\end{notation}

The following lemma is similar in essence to \cite[Lemma~1]{Libk95}:

\begin{lemma}\label{L:N(L)commutes}
Let $L$ be a bounded lattice, let $a\in\Neu L$, let $x$, $y\in L$. Then
$y=a\vee x$ \pup{resp., $y=a\wedge x$} implies that $y=a\vee^*x$
\pup{resp., $y=a\wedge^*x$}.
\end{lemma}

\begin{proof}
We prove, for example, that $y=a\vee x$ implies that $y=a\vee^*x$. Let
$p\in\J(L)$ such that $p\leq y$. Then, by using the fact that $a$ is neutral,
$p=p\wedge(a\vee x)=(p\wedge a)\vee(p\wedge x)$, hence, since $p$ is \jirr,
either $p\leq a$ or $p\leq x$. The proof for the meet is similar.
\end{proof}

We leave to the reader the straightforward proof of the following lemma:

\begin{lemma}\label{L:sqtononsq}
Let $L$ be a \fbs\ bounded lattice.
Let $x$, $y\in L$, let $(x_i)_{i\in I}$ be a family of elements of $L$. Then
the following assertions hold:
\begin{enumerate}
\item $x=\bigvee^*_{i\in I}x_i$ implies that $x=\bigvee_{i\in I}x_i$;

\item $x=\bigwedge^*_{i\in I}x_i$ implies that $x=\bigwedge_{i\in I}x_i$;

\item $x=\bigvee^*_{i\in I}x_i$ implies that
$x\wedge y=\bigvee^*_{i\in I}(x_i\wedge y)$;

\item $x=\bigwedge^*_{i\in I}x_i$ implies that
$x\vee y=\bigwedge^*_{i\in I}(x_i\vee y)$.
\end{enumerate}
\end{lemma}

\begin{lemma}\label{L:CplPair}
Let $L$ be a \fbs\ bounded lattice.
Let $a$, $b\in L$. Then the following are equivalent:
\begin{enumerate}
\item $a\vee^*b=1$ and $a\wedge^*b=0$;

\item $(a,b)$ is a complementary pair of elements of $\Cen L$.
\end{enumerate}

\end{lemma}

\begin{proof}
(i)$\Rightarrow$(ii) We consider the maps $f\colon L\to[0,a]\times[0,b]$ and
$g\colon[0,a]\times[0,b]\to L$ defined by the following formulas:
 \begin{align*}
 f(z)&=(z\wedge a,z\wedge b),&&\text{for all }z\in L,\\
 g(x,y)&=x\vee y,&&\text{for all }(x,y)\in[0,a]\times[0,b].
 \end{align*}
For any $z\in L$, it follows from Lemma~\ref{L:sqtononsq} that
$z=(z\wedge a)\vee(z\wedge b)$, so
$g\circ f=\id_L$. Conversely, let $x\leq a$ and $y\leq b$ in $L$. Then, again
by using Lemma~\ref{L:sqtononsq},
$x\leq(x\vee y)\wedge a\leq(x\vee b)\wedge(x\vee a)=x\vee(a\wedge^*b)=x$,
whence $x=(x\vee y)\wedge a$. Similarly, $y=(x\vee y)\wedge b$. Therefore,
$f\circ g=\id_{[0,a]\times[0,b]}$, so $f$ and $g$ are mutually
inverse isomorphisms. The conclusion (ii) follows then from
Lemma~\ref{L:charactC(L)}.

(ii)$\Rightarrow$(i) follows immediately from Lemma~\ref{L:N(L)commutes}.
\end{proof}

Now we can prove one of the main lemmas of this section:

\begin{lemma}\label{L:CplfbsC(L)}
Let $L$ be a \fbs\ complete lattice. Then the center $\Cen L$ is a complete
sublattice of $L$.
\end{lemma}

\begin{proof}
Let $(a_i)_{i\in I}$ be a family of elements of $\Cen L$, then
put $a=\bigvee_{i\in I}a_i$ and $b=\bigwedge_{i\in I}\neg a_i$.

We first claim that $a\vee^*b=1$ and $a\wedge^*b=0$. Indeed, let us prove
for example the first assertion. Let $p\in\J(L)$. It follows from 
Lemma~\ref{L:CplPair} that for all $i\in I$, either $p\leq a_i$ or
$p\leq\neg a_i$. Hence, if $p\nleq b$, then there exists $i\in I$ such that
$p\nleq\neg a_i$, whence $p\leq a_i\leq a$. The proof of $a\wedge^*b=0$ is
dual. It follows, again by Lemma~\ref{L:CplPair}, that $(a,b)$ is a
complementary pair of $\Cen L$. In particular, $\Cen L$ is a complete
sublattice of $L$.
\end{proof}

By Lemma~\ref{L:CplfbsC(L)}, for any $x\in L$, there is a least element
$u$ of $\Cen L$ such that $x\leq u$, we denote this element by $e(x)$, the
\emph{central cover} of $x$.

\begin{lemma}\label{L:C(L)atomic}
Let $L$ be a \fbs\ complete lattice.
The Boolean lattice $\Cen L$ is atomistic, with atoms the $e(p)$ for
$p\nobreak\in\nobreak\J(L)$.
\end{lemma}

\begin{proof}
Observe first the obvious equality
$u=\bigvee\setm{e(p)}{p\in\J(L),\ p\leq u}$,
for any $u\in\Cen L$. Hence, it suffices to prove that $e(p)$ is
an atom of $\Cen L$, for all $p\in\J(L)$. Suppose otherwise. Then
$e(p)=u\vee v$, for nonzero elements $u$ and $v$ of $\Cen L$ such that
$u\wedge v=0$. {}From $u\in\Neu L$ follows that
$p=p\wedge(u\vee v)=(p\wedge u)\vee(p\wedge v)$, whence, since $p\in\J(L)$,
either $p\leq u$ or $p\leq v$. Suppose, for example, that $p\leq u$. Then
$e(p)\leq u$, whence $v=0$, a contradiction.
\end{proof}

{}From Proposition~\ref{P:SuffDec}(iii) and Lemmas \ref{L:CplfbsC(L)} and
\ref{L:C(L)atomic}, we can now deduce immediately the main result of this
section:

\begin{theorem}\label{T:DecompL}
Every \fbs\ complete lattice is isomorphic to a direct product of directly
indecomposable lattices.
\end{theorem}

As immediate corollaries of Theorem~\ref{T:DecompL} and the fact that every
algebraic lattice is dually spatial (see \cite[Theorem I.4.22]{Comp}, or
\cite[Lemma~1.3.2]{Gorb}), we observe the following, see
\cite[Theorem~2]{Libk95}:

\begin{corollary}[Libkin's Decomposition Theorem]
\label{C:Libkin}
Every algebraic and spatial lattice is isomorphic to a direct product of
directly indecomposable lattices.
\end{corollary}

In particular, every algebraic and atomistic lattice is isomorphic to a
direct product of directly indecomposable lattices. In fact, since every
algebraic lattice is dually spatial, Theorem~\ref{T:DecompL} makes it
possible to extend Corollary~\ref{C:Libkin} to finitely spatial algebraic
lattices. In particular, we obtain the following consequence, a stronger
form of which is stated in \cite[Corollary~2]{Wale97}:

\begin{corollary}\label{C:AtduAt}
Every algebraic and dually algebraic lattice is isomorphic to a direct
product of directly indecomposable lattices.
\end{corollary}

\begin{example}\label{Ex:DCX}
There exists a dually algebraic, atomistic, distributive lattice $D$ whose
center $\Cen D$ is not complete. In particular, $D$ cannot be decomposed as
a direct product of directly indecomposable lattices.
\end{example}

\begin{proof}
We recall that the \emph{interval topology} on a \emph{totally} ordered set
$T$ is the least topology on $T$ for which all intervals of the form $[a)$
(resp., $(a]$) are closed subsets. It is a well-known result, due to
O.~Frink (see for example \cite[Theorem~X.12.20]{Birk}), that states
that the interval topology on $T$ is compact Hausdorff if{f} $T$ is a
complete lattice.

Now we endow the ordinal $\omega+1$ with its interval topology, and we let
$D$ be the lattice of all closed sets of this topology. Hence,
 \[
 D=\setm{x\subset\omega}{x\text{ is finite}}\cup
 \setm{x\cup\set{\omega}}{x\subseteq\omega}.
 \]
Observe that $D$ is a closure system in the powerset algebra $\Pow(\omega+1)$
of $\omega+1$, thus it is a complete lattice. Moreover, $D$ is a distributive
sublattice of $\Pow(\omega+1)$, and it is atomistic since every element of $D$
is a union of singletons. Moreover, it is straightforward to compute that
 \[
 \Cen D=\setm{x\subset\omega}{x\text{ is finite}}\cup
 \setm{x\cup\set{\omega}}{x\subseteq\omega\text{ is cofinite}},
 \]
so $\Cen D$ consists exactly of the \emph{clopen} subsets of $\omega+1$.
Since $\omega+1$ is a compact topological space, every element of $\Cen D$
is dually compact in $D$. Furthermore, every closed subset of $\omega+1$
is an intersection of clopen subsets, therefore, \emph{$D$ is dually
algebraic}.

Put $a=\setm{2m+1}{m<\omega}$ and $b_n=(\omega+1)\setminus\set{2n}$, for
all $n<\omega$. Observe that both $a\cap m$ and $b_n$ belong to $\Cen D$,
for all $m$, $n<\omega$, and that $a\cap m\subset b_n$. However, there is
no element $x$ of $\Cen D$ such that $a\cap m\subseteq x\subseteq b_n$ for
all $m$, $n<\omega$, because otherwise either $x=a$ or
$x=a\cup\set{\omega}$ would belong to $\Cen D$, a contradiction.
\end{proof}

\begin{remark}
It is easy to verify that $D$ is even \emph{strongly atomic}, that is,
$a<b$ implies that there exists $x\in D$ such that $a\prec x\leq b$, for all
$a$, $b\in D$. We recall that every algebraic
lattice $A$ is \emph{weakly atomic}, that is, for all $a<b$ in $A$, there
are $x$, $y\in A$ such that $a\leq x\prec y\leq b$ (see
\cite[Lemma~2.2]{CrDi} or \cite[Exercise~1.3.1]{Gorb}).
\end{remark}

\section{Direct decompositions of lattices of algebraic subsets}
\label{S:AlgSubs}

For a complete lattice $A$, a subset $X$ of $A$ is \emph{algebraic}, if $X$
is closed under arbitrary intersections and nonempty up-directed joins, and we
denote by $\SP(A)$ the lattice of all algebraic subsets of $A$. Then the
following basic lemma holds, see~\cite{AGT} for more information:

\begin{lemma}
Let $A$ be a complete lattice. Then the following assertions hold:
\begin{enumerate}
\item If $A$ is upper continuous, then $\SP(A)$ is a \jsd, lower continuous
lattice.

\item If $A$ is algebraic, then $\SP(A)$ is dually algebraic.

\item If $A=\Pow(X)$ for some set $X$, then $\SP(A)$ is dually spatial.
\end{enumerate}
\end{lemma}

We recall at this point that any algebraic lattice is upper continuous (see
\cite[Lemma~2.3]{CrDi}), and that for $A$ a general algebraic lattice,
$\SP(A)$ does not need to be dually spatial (see \cite{AGT}). We also observe
that $\SP(A)$ is the closure lattice of the atomistic closure space
$(A\setminus\set{1},\SPl)$, where, for every subset $X$ of $A$,
we put $\SPl(X)=\ol{X}\setminus\set{1}$, where $\ol{X}$ denotes the
algebraic subset of $A$ generated by $X$.

K.\,V. Adaricheva has kindly informed the author that all lattices of the
form $\SP(\Pow(X))$ are directly indecomposable. A stronger result is the
following:

\begin{proposition}\label{P:SPIndec}
For any complete Boolean algebra $B$, the lattice $\SP(B)$ is
subdirectly irreducible.
\end{proposition}

\begin{proof}
The atoms of $\SP(B)$ are the $U_a=\set{a,1}$ for
$a\in B\setminus\set{1}$. For $X$, $Y\in\SP(B)$,
we denote by $\Theta(X,Y)$ the principal congruence of
$\SP(B)$ generated by the pair $(X,Y)$.
For any $a\in B\setminus\set{0,1}$, the containment
$U_0\subset U_a\vee U_{\neg a}$ holds, with $U_0$, $U_a$, and $U_{\neg a}$
distinct atoms of $\SP(B)$, thus
$\Theta(\set{1},U_0)\subseteq\Theta(\set{1},U_a)$.
Since $\SP(B)$ is atomistic, it follows that $\Theta(\set{1},U_0)$ is the
smallest nonzero congruence of~$\SP(B)$.
\end{proof}

\begin{remark}\label{Rk:Notsimple}
The lattice $\SP(\Pow(2))$ is the (finite) \jz-semilattice defined by
generators $a$, $b$, $c$ and the unique relation $c\leq a\vee b$, hence it is
not simple.
\end{remark}

\begin{remark}\label{Rk:NotSI}
Even for finite atomistic lattices which are lower bounded homomorphic
images of free lattices, direct indecomposability is not equivalent to
subdirect irreducibility. For example, the lattice $\Co(\two^2)$ of all
convex subsets of $\two^2$ (diagrammed for example in \cite[p.~224]{BB})
is directly indecomposable, although not subdirectly irreducible. This is
another strong point of contrast between geometric lattices and convex
geometries (see \cite{AGT} for the latter): namely, every directly
indecomposable geometric lattice is subdirectly irreducible, see
\cite[Theorem~IV.3.6]{GLT2}.
\end{remark}

On the other hand, as we shall see in a moment, the lattice $\SP(A)$ displays
a very different behavior for $A$ a \emph{totally ordered} algebraic lattice.
The proof of the following lemma is a straightforward exercise.

\begin{lemma}\label{L:TopAlgClos}
Let $A$ be a totally ordered algebraic lattice. Then a subset $X$ of $A$
belongs to $\SP(A)$ if{f} $X$ is closed for the interval topology and
$1\in X$.
\end{lemma}

\begin{example}[see \cite{ADG2}]\label{Ex:NonCpleteCen}
Let $C=[0,1]$ be the rational unit interval, let $A$ be the ideal lattice of
$C$. Then $\Cen\SP(A)$ is not a complete lattice.
\end{example}

\begin{proof}
Put $j(x)=[0,x]$, for all $x\in C$, and, if $x>0$, put $j(x)_*=[0,x)$,
so $j(x)_*\prec j(x)$.
Observe that $A=\setm{j(x)}{x\in C}\cup\setm{j(x)_*}{x\in C\setminus\set{0}}$
is a complete chain with top element $\one=[0,1]$. It follows from
Lemma~\ref{L:TopAlgClos} that $\SP(A)$ is isomorphic to the lattice $D$ of
all closed subsets of $A\setminus\set{\one}$ endowed with the interval
topology. Therefore, the center $\Cen\SP(A)$ is isomorphic to the Boolean
lattice $B$ of all clopen subsets of $A\setminus\set{\one}$ for the interval
topology.

Now put $a_n=\frac{1}{2}-\frac{1}{2n}$ for every positive integer $n$, and
$a=\frac{1}{2}$. For each positive integer $n$, we put
 \[
 X_n=[j(a_{2n}),j(a_{2n+1})_*],\qquad
 Y_n=\bigcup_{0<k\leq n}X_k\cup[j(a_{2n+2}),j(a)_*].
 \]
Then both $X_m$ and $Y_n$ are clopen subsets of
$A\setminus\set{\one}$ with $X_m\subset Y_n$, for all $m$, $n>0$. However,
the only subsets $Y$ of $A\setminus\set{\one}$ such that
$X_m\subseteq Y\subseteq Y_n$ for all $m$, $n>0$ are
$Z=\bigcup_{0<k<\omega}X_k$, which is not closed, and $Z\cup\set{j(a)_*}$,
which is not open.
\end{proof}

We conclude the paper with a problem:

\begin{problem}
Find a common generalization of Theorem~\ref{T:DecompL} (decomposition
theorem for finitely bi-spatial complete lattices) and various
decomposition results such as the ones in \cite{Jaku96,Jano67,Wale97}.
\end{problem}

Indeed, the hard core of Theorem~\ref{T:DecompL} and its analogues lies in
proving that the center is complete. All the methods used here and in
\cite{Jaku96,Jano67,Wale97} bear some formal similarity, but none of the
results seems to follow from the others.

\section*{Acknowledgment}
The author is grateful to Kira Adaricheva for having gotten him
interested in the topic and for her many helpful comments.

\end{document}